\documentclass[12pt]{article}
\usepackage{makeidx,epsfig,setspace,graphicx,amsfonts,amsmath,amssymb,
amsthm,latexsym,citesort}
\usepackage{amsmath}
\usepackage{graphicx}
\usepackage{amsfonts}
\usepackage{amssymb}
\makeindex
\newtheorem{theorem}{Theorem}

\newtheorem{remark}[theorem]{Remark}

\newcommand{\RR}{\mathbb{R}}

\begin{document}

\title{On reconstruction
formulas and algorithms for the thermoacoustic tomography}
\date{}

%\chapter[Reconstruction formulas for TAT]{On reconstruction
%formulas and algorithms for the thermoacoustic tomography}

%\author{M.~Agranovsky, P.~Kuchment, and L.~Kunyansky}

\author{
Mark Agranovsky\\
Department of Mathematics\\
Bar-Ilan University, Israel\\
\\
Peter Kuchment\\
Department of Mathematics\\
Texas A\& M University, USA\\
\\
Leonid Kunyansky\\
Department of Mathematics\\
University of Arizona, USA}

\maketitle

%\pagenumbering{roman}
%\begin{doublespace}
%\frontmatter
\tableofcontents
%\mainmatter
%\pagenumbering{arabic}

\section*{Introduction}
%%%%%%%%%%%%%%%%%%%%%%%
Recent years have brought about exciting new developments in
computerized tomography. In particular, a novel, very promising
approach to the creation of diagnostic techniques
consists in combining different imaging modalities, in order
to take advantage of their individual strengths.
Perhaps, the most successful example of such a combination is
the \textbf {Thermoacoustic Tomography (TAT)}
(also called photoacoustic tomography and optoacoustic tomography and
abbreviated as TCT, PAT, or
OAT)~\cite{Kruger,Kruger03,Oraev96,Oraev,Oraev2,PAT,Wang_book,MXW_review}.

Major progress has been made recently in developing
the mathematical foundations of TAT, including proving
uniqueness of reconstruction, obtaining range descriptions for
the relevant operators, deriving inversion formulas and
algorithms, understanding solutions of incomplete data problems,
stability of solutions,
etc. One can find a survey of these results and extensive
bibliography in~\cite{KuKu}. In the present article we
concentrate on the recent advances in the inversion
formulas and algorithms for TAT. Mathematical problems of the
same type arise also in sonar, radar, and geophysics applications
(e.g.,~\cite{LQ,NC,Beylkin}). Discussion of some mathematical problems
concerning TAT can be also found in the chapters written by
D.~Finch and Rakesh and by S.~Patch.

While this text addresses the mathematics of TAT only, one can find
extensive discussion of physics, engineering, and biological issues related to
TAT in the recent surveys \cite{Oraev,Oraev2,MXW_review}, textbook \cite{Wang_book},
as well as in other chapters of this volume.

%%%%%%%%%%%%%%%%%%%%%
\section{Thermoacoustic tomography}
%%%%%%%%%%%%%%%%%%%%%

We give first a brief description of TAT. The data acquisition starts with
a short electromagnetic (EM) pulse being sent through the biological object
under investigation (e.g., woman's breast in mammography)
\footnote{It has been argued that the radiofrequency and visible light
ranges are most appropriate in TAT \cite{MXW_review}. For the purpose of
this text, no distinction is made between these cases.}.
\begin{figure}[th]
\begin{center}
%\scalebox{0.7}
{\includegraphics{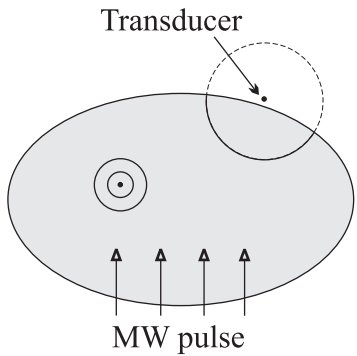}}
\end{center}
\caption{The TAT procedure.}%
\label{F:tat}%
\end{figure}
\noindent
%It is assumed that the whole object is irradiated with sufficient
%uniformity.
A fraction
%$f(x)$
of EM energy is absorbed at each
location $x$ inside the object, thus triggering thermoelastic
expansion of the tissue and emergence of a pressure
wave $p(x,t)$ (an ultrasound signal) that, in turn, is measured
by transducers placed along some \textit{observation surface}  $S$
surrounding (completely or partially) the object.
The initial pressure $p_0(x)= p(x,0)$ is determined by the
intensity of the EM pulse (that assumed to be known) and by the
local properties of the tissue. It is known
(e.g.,~\cite{Kruger,Oraev,Oraev2,MXW1,MXW_review}) that in the
radiofrequency and visible light ranges
absorption of the EM energy by cancerous cells
is several times stronger than by the healthy ones.
Thus, knowledge of the initial pressure
%function $f(x)$ describing
%spatial distribution of the absorbed energy
$p_0(x)$
would provide
an efficient tool for early detection of cancer.
Frequently, the ultrasound contrast is sufficiently small
to justify the use of the constant sound speed approximation.
Most work on TAT up to date is based on this assumption.
However, such an approximation is not always appropriate; some
of the results described below, as well as in~\cite{AKu,JinWang,KuKu}
aim towards the general case of a variable speed
of sound.

Once the data $p(x,t)$ has been measured on $S\times\RR^+$,
one can attempt to recover from $p(x,t)$ the initial value $p_0(x)$ of
the pressure inside $S$ (the thermoacoustic image).

\section[Mathematical model of TAT]{Mathematical model of TAT}
%%%%%%%%%%%%%%%%%%%%%
Let us for notational convenience denote $p_0(x)$
(the image to be reconstructed) by $f(x)$.
In this section, we present a mathematical description of the
relation between the functions $f(x)$  and $p(x,t)$.
We assume that the function $f(x)$
is compactly supported in $\RR^n$ (we allow the dimension to be
arbitrary, albeit the most interesting cases for TAT are $n=3$
and $n=2$). At each point $y$ of an {\bf observation
surface} $S$ one places a point detector\footnote{Planar and linear
detectors have been considered as well, see
\cite{haltmaier_large,haltmaier_fabri} and further references
in~\cite{KuKu}.} that measures the value of the pressure $p(y,t)$ at any
moment $t>0$. It is usually assumed that the surface $S$ is closed (rather
than, say, cylinder or a plane\footnote{Reconstruction formulas for the
planar and cylindrical cases are well known, see e.g.
\cite{Natt_new,And,Faw,XFW,XXW}.}). It is also assumed that the object
(and thus the support of $f(x)$) is completely surrounded by $S$.
The latter assumption is crucial for the validity of most
inversion formulas; however in some cases we will be able
to abandon this requirement.

The mathematical model described below relies upon some physical
assumptions on the measurement process, which we will not describe here.
The reader can find such a discussion in \cite{MXW_review}.

We assume that the ultrasound speed $v_s(x)$ is known, e.g.,
through transmission ultrasound measurements~\cite{JinWang}.
Then, the pressure wave $p(x,t)$ satisfies the following set
of equations~\cite{Diebold,Tam,MXW1}:
\begin{equation}\label{E:wave}
\begin{cases}
    p_{tt}={v_s}^2(x)\Delta_x p, \quad t\geq 0, \quad x\in\RR^n\\
    p(x,0)=f(x),\\
    p_t(x,0)=0
    \end{cases}
\end{equation}
Now one needs to recover the initial value $f(x)$ at $t=0$ of the
solution $p(x,t)$ from the measured data $g(y,t):=p(y,t),  y\in S,  t\geq 0$.
Incorporating this data, one rewrites~(\ref{E:wave}) as
\begin{equation}\label{E:wave_data}
\begin{cases}
    p_{tt}={v_s}^2(x) \Delta_x p, \quad  t\geq 0, \quad x\in\RR^n\\
    p(x,0)=f(x),\\
    p_t(x,0)=0\\
    p(y,t)=g(y,t), \quad y\in S\times\RR^+
    \end{cases}
\end{equation}

\begin{figure}[ht]
\begin{center}
%\scalebox{0.7}
{\includegraphics{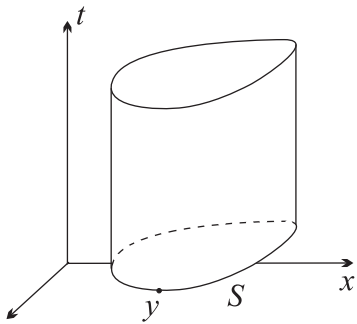}}
\label{F:cylinder}
\end{center}
\caption{An illustration to~(\ref{E:wave_data}).}
\label{F:wave}
\end{figure}
In other words, we would like to recover the initial value $f(x)
$ in~(\ref{E:wave_data}) from the knowledge of the lateral data
$g(y,t)$ (see Figure~\ref{F:cylinder}). At a first glance, it
seems that the data is insufficient for the reconstruction, i.e. for
recovering the solution of the wave equation in a cylinder from
the lateral values alone. However, this impression is incorrect,
since there is additional information
that the solution holds in the whole space, not just inside the
cylinder $S\times \RR^+$. To put it differently, if one solves not
only the internal, but also the external problem for the wave equation
with the data $g$ on the cylinder $S\times\RR^+$, then the solutions must
have matching normal derivatives on $S\times\RR^+$. In most cases, this
additional information provides uniqueness of recovery of $f(x)$ (see
below, as well as
%~\cite{ABK}-\cite{AmbKuc_inj}, ~\cite{KuKu,AKu,FPR}
~\cite{ABK,AKQ,AQ,AmbKuc_inj,KuKu,AKu,FPR}, and
references therein). It is also sometimes useful to notice
that $p$ can be extended as an even function of
time and thus satisfies the wave equation for all values of $t$. Similarly,
data $g$ can be extended to an even function. This, in particular enables
one to apply Fourier transform in time.

An additional structure arises in this problem, if one
assumes that the object under investigation is nearly
homogeneous with respect to ultrasound: $v_s(x)=1$.
In this constant speed case, there is an alternative way to
describe the relation between the data $g(y,t),
(y,t)\in S\times\mathbb{R}^{+}$ and the unknown image $f(x),
x\in\mathbb {R}^{3}$. The known
Poisson-Kirchhoff formulas~\cite[Ch. VI, Section 13.2, Formula
(15)]{CH} for the solution of~(\ref{E:wave}) with $v_s=1$ give%
\begin{equation}
p(x,t) = \frac{\partial}{\partial t}\left(  t(Rf)(x,t)\right)  ,
\label{E:KP}%
\end{equation}
where
\begin{equation}
(Rf)(x,r)= \frac{1}{4\pi}\int\limits_{|y|=1} f(x+ry)dA(y) \label{E:mean}%
\end{equation}
is the \emph{spherical mean operator} applied to the function $f(x)$,
and $dA$ is the surface area element on the unit sphere in $\mathbb
{R}^{3}$.
Thus, the function $g(y,t)$ for $y\in S$ and all $t\geq0$
essentially carries the same information as the
spherical mean $Rf (y,t)$ at all points
$(y,t)\in S\times\mathbb{R}^{+}$ (see, e.g., \cite{AQ}). One can, therefore, study the
spherical mean operator $R:f\to Rf$ and, in particular, its
restriction $R_{S}$ to the
points $y\in S$ of the observation surface:
\begin{equation}
R_{S}f(x,t)=\int\limits_{|y|=1}f(x+ty)dA(y), \quad x\in S, \quad t\geq0. \label{E:Radon_S}%
\end{equation}
This explains why in many studies on thermoacoustic tomography,
the spherical mean operator has been used as the model.
One needs to notice, though, that in the case of a non-constant
sound speed, the spherical mean interpretation (as well as any
integral geometry approximation) is no longer valid,
while the wave equation model still is.

%%%%%%%%%%%%%%%%%%%%%
\section{Uniqueness of reconstruction}
%%%%%%%%%%%%%%%%%%%%%

Uniqueness of reconstruction of a compactly supported (or
sufficiently  fast decaying) function $f(x)$ from the data $g$
collected from a closed surface $S$ is well known in the case of
a constant sound speed (i.e., when the interpretation in terms of
spherical mean operators is possible). One can find discussion of
such results in~\cite{ABK,AKu,AQ,AmbKuc_inj,FPR,Kuc93,Kuch_AMS05,KuKu,FR2,FR3}.

In the case of a variable sound speed, it is shown in \cite[Theorem 4]{FR2}
that uniqueness of reconstruction also holds for a smoothly varying
(positive) sound speed, if the function $f(x)$ is supported inside the
observation surface $S$. The proof uses the famous unique continuation
theorem by D.~Tataru \cite{Tataru}.

We present now a recent simple uniqueness theorem that also allows a
non-constant sound speed $v_s(x)$ and does not require the function to be
supported inside $S$. In order to do so, we need to formulate first some
assumptions on $v_s(x)$ and the function $f(x)$ to be reconstructed.

\begin{enumerate}
\item Support of $f(x)\in H^s_{loc}(\RR^n), s>1/2$ is compact.

\item The sound speed is smooth (a condition that can be reduced),
strictly positive $v_s(x)>v_0>0$ and such that $v_s(x)-1$ has
compact support, i.e. $v_s(x)=1$ for large $x$.

\item Consider the Hamiltonian system in $\RR^{2n}_{x,\xi}$ with the
Hamiltonian $H=\frac{{v_s}^2(x)}{2}|\xi|^2$:
\begin{equation}
\begin{cases}
x^\prime_t=\frac{\partial H}{\partial \xi}={v_s}^2(x)\xi \\
\xi^\prime_t=-\frac{\partial H}{\partial x}=-\frac 12 \nabla
\left({v_s}^2(x)\right)|\xi|
^2 \\
x|_{t=0}=x_0, \quad \xi|_{t=0}=\xi_0.
\end{cases}
\end{equation}
The solutions of this system are called {\em bicharacteristics} and
their projections into $\RR^n_x$ are {\em rays}.

We will assume that the {\bf non-trapping condition} holds, i.e.
that  all rays (with $\xi_0\neq 0$) tend to infinity when $t \to
\infty$.

\end{enumerate}

\begin{theorem}\cite{AKu}\label{T:uniqueness}
Under the assumptions formulated above, compactly
supported function $f(x)$ is uniquely determined by the data $g$.
(No assumption of $f$ being supported inside $S$ is imposed.)
\end{theorem}

Uniqueness fails, however, if $f$ does not decay sufficiently
fast (see~\cite{ABK}, where it is shown for the constant speed
in which spaces $L^{p}(\mathbb{R}^{d})$ of
functions $f(x)$ closed surfaces remain uniqueness sets).

%%%%%%%%%%%%%%%%%%%%%
\section[Reconstruction:constant speed]{Reconstruction in the case of constant sound speed: formulas,
algorithms, and examples.}
%%%%%%%%%%%%%%%%%%%%%
%
%
%
We consider here the case of a constant sound speed: $v_s(x)=1$.
One can work then either with the wave equation, or with the
spherical mean operator model.

\subsection{Inversion formulas and procedures}

Consider the case of the observation surface $S$ being a sphere.
The first inversion procedures for this situation were
obtained in~\cite{Norton1} in $2D$ and in~\cite{Norton2} in $3D$
by harmonic decomposition of the measured data $g$ and of
the function $f$, and then by equating coefficients of the
corresponding Fourier series (see also \cite{KuKu} for a brief description
of this procedure). The two resulting series solutions are not quite analogous.
Indeed, in~\cite{Norton1} one had to divide the Hankel transform of the
data by the Bessel functions that have infinitely many zeros, which would
create instabilities during implementation. The $3D$ solution in~\cite{Norton2}
is free of this difficulty and can also be adopted for $2D$. We will see a
different type of series solutions later on in this section.

\subsubsection{Approximate inversion formulas}

The standard way of inverting Radon transform in tomographic
applications is by using filtered backprojection type
formulas~\cite{Leon_Radon,GGG,Helg_Radon,Natt_old,Natt_new}.
It combines
a linear filtration of projections (either in Fourier domain, or by
a convolution with a certain kernel) followed (or preceded) by a
backprojection. In the case of the set of spheres centered on a
closed surface (e.g., sphere) $S$, one expects such a formula to
involve a filtration with respect to the radial variable and an
integration over the set of spheres passing through the point $x$ of
interest. Albeit for quite a long time no such formula had been
discovered, this did not prevent practitioners from reconstructions.
The reason was that good approximate inversion formulas
(parametrices) could be developed, followed by an optional iterative
improvement of the reconstruction
\cite{PopSush,PopSush2,PAT,MXW1,XFW,XXW,XWAK}.

Perhaps the most advanced approach of this kind was adopted by
Popov and Sushko~\cite{PopSush,PopSush2}. These authors have
developed a set of
''straightening'' formulas that allow one to reconstruct from
the spherical
means an approximation to the regular Radon projections. The
main idea is that for each (hyper)plane passing through the
support of the function to be
reconstructed, one builds a family of spheres with centers at
the detectors'
locations and tangential to that plane. One such sphere is
chosen for each
point of the plane contained within the support. The integrals
over these
spheres are known, as they form a subset of projections $g$. An
approximation to the integral of the function over the plane is
then computed by integrating over these projections a functional
(local in odd and non-local in even dimensions). When all the
plane integrals are computed,
the function is reconstructed by applying inversion formulas for
the regular Radon transform. This procedure is not exact;
however, as shown in~\cite{PopSush},
such an algorithm yields a parametrix. Namely,
the difference between such an approximation and the original
function $f$ is
described by a pseudodifferential operator of order $-1$ applied
to $f$. In other words,
reconstruction is accurate up to a smoothing operator. This
result holds even
if the measuring surface is not closed (but satisfies a
''visibility''
condition), which is important for applications in the problems
with incomplete data.

\subsubsection{Exact filtered backprojection formulas in $3D$}

The first set of exact inversion formulas of the filtered
backprojection type for the spherical surface $S$ was discovered
in~\cite{FPR}. These formulas were obtained only in odd dimensions
(and then extended to even dimensions in~\cite {Finch_even}).
Various versions of such formulas (different in terms of the
order in which the filtration and backprojection steps are
performed) were developed.

To describe these formulas, let us assume that $B$ is the unit
ball, $S=\partial B$ is the unit sphere in $\mathbb{R}^{3}$, and
a function $f(x)$ is supported inside $S$. The values of its
spherical integrals $g(z,r)$ with the centers on $S$ are assumed
to be known:
\begin{equation}
g(z,r)=\int\limits_{\mathbb{S}^{2}}f(z+rs)r^{2}dA(s)=4\pi r^{2}%
R_{S}f(z,r),\qquad z\in S.\label{E:spher_int}%
\end{equation}
Some of the $3D$ inversion formulas of~\cite{FPR} are:
\begin{align}
f(y) &  =-\frac{1}{8\pi^{2}
%R
}\Delta_{y}\int\limits_{S}%
\frac{
g(z,|z-y|)}{|z-y|}dA(z),\label{FPR3da}\\
f(y) &  =-\frac{1}{8\pi^{2}
%R
}\int\limits_{S}\left
(  \frac{1}{t}%
\frac{d^{2}}{dt^{2}}g(z,t)\right)  \left.  {\phantom{\rule{1pt}
{8mm}}}\right|
_{t=|z-y|}dA(z).\label{FPR3d}%
\end{align}

A different set of explicit inversion formulas, which work in
arbitrary dimensions, was found in~\cite{Kunyansky}. In 3D case
the general
expression derived in~\cite{Kunyansky} simplifies to
\begin{equation}
f(y)=\frac{1}{8\pi^{2}}\mathrm{div}\int\limits_{S}n%
(z)\left(  \frac{1}{t} \frac{d}{dt}\frac{g(z,t)}{t}\right)  \left.  {\phantom
{\rule
{1pt}{8mm}}}\right|  _{t=|z-y|}dA(z), \label{E:kunyansky}%
\end{equation}
where $n(z)$ is the vector of exterior normal to $S$.  (We eliminated
in this expression the minus sign erroneously present in the original
formula.)
Equation~(\ref{E:kunyansky}) is equivalent to one of the $3D$ formulas
derived earlier in~\cite{MXW2}.

Similarly to the case of the standard ``flat'' Radon transform,
all these $3D$ inversion formulas are
local, i.e. in order to reconstruct a value of the function at a
certain point,
one needs to know only values of all the integrals over the
spheres passing
through an infinitesimally small neighborhood of that point.

It is worth noting that although formulas~(\ref{FPR3d}) and~(\ref
{E:kunyansky}) yield identical results when applied to functions
that belongs to the
range of the spherical mean Radon transform, they are in general
not equivalent, i.e. lead to different reconstructions when the
data is outside of the range (for instance, due to errors).
Another important fact about these reconstruction techniques is
that, unfortunately, they do not yield correct
reconstruction within the region surrounded by the detectors if
the source is not contained within this region. Both these
statements can be easily proven by the following example. Let us
assume that the source function $f(x)$ is constant
(equal to 1) within the ball $B(0,3)$ of radius 3 centered at
the origin. In order to reconstruct the function within the unit
ball, both formulas~(\ref{E:kunyansky}) and~(\ref{FPR3d}) use
only integrals over spheres with the radius less or equal to 2,
and centered at the points of the unit sphere.
Obviously, all these spheres lie within the $B(0,3)$, and thus
the projections $g(z,t)$ are equal to the areas of the
corresponding integration spheres, i.e. to $4\pi t^{2}.$ By
substituting this expression into~(\ref{FPR3d}), we obtain
\[
f_{1}(y)=-\frac{1}{\pi R}\int\limits_{S}%
\frac{1}{|z-y|}dA(z).
\]
Function $f_{1}(y)$ defined by the above formula is harmonic
in the interior of $B$, since the integrand is
the free space Green's function of the Laplace equation. Due to
the symmetry of the
geometry, $f_{1}(y)$ is radially symmetric, i.e. it depends only
on $|y|$. Therefore $f_{1}(y)=const$ for all $y\in B\setminus
S$. Let us compute $f_{1}(0)$:
\[
f_{1}(0)=-\frac{1}{\pi R}\int\limits_{S}\frac{1}{R}
dA(z)=-4.
\]
Thus, $f_{1}(y)=-4$ for all $y\in B\setminus S$.

A similar computation with the use of~(\ref{E:kunyansky}) yields
\begin{align*}
f_{2}(y)  &  = \frac{1}{2\pi}\mathrm{div}\int
\limits_{S}n(z)\frac{1}{|z-y|}dA(z)\\
&  =-\frac{1}{2\pi}\int\limits_{S}\frac{d}{dn(z)}\frac
{1}{|z-y|}dA(z) =\frac{4\pi}{2\pi}=2,
\end{align*}
where we used the $3D$ Gauss formula. Both results $f_1$ and $f_2$
are incorrect (not equal to $1$). Besides, they are different,
which proves that formulas~(\ref{FPR3d}) and~(\ref{E:kunyansky})
are not equivalent.

One of the important benefits of having exact inversion formulas
is that often a rather straightforward discretization of such a formula
yields an efficient and stable
reconstruction algorithm. Such algorithms were developed in~\cite
{AmbPatch} using equations~(\ref{FPR3da}) and~(\ref{FPR3d}),
and in~\cite{Kunyansky} utilizing formula~(\ref
{E:kunyansky}).

\begin{figure}[th]
\begin{center}
%\scalebox{0.7}
%
%
%
{\includegraphics{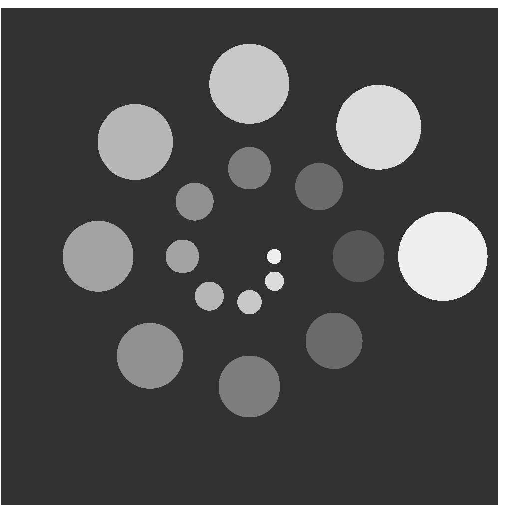} \phantom{aaa} \includegraphics
{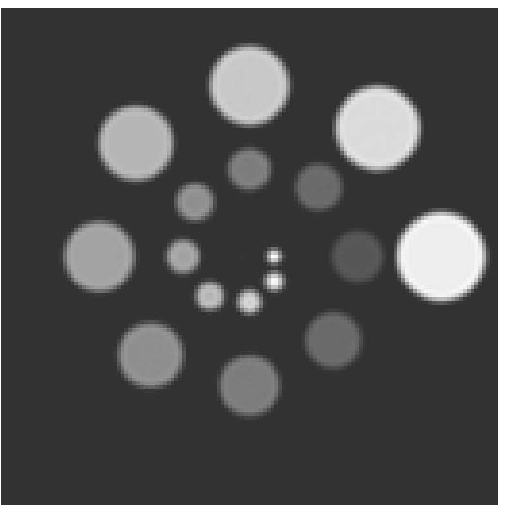}}
\end{center}
\caption{A mathematical phantom in $3D$ (left) and its reconstruction using
inversion formula~(\ref{E:kunyansky})}%
\label{F:recon}%
\end{figure}

In the simplest case, when the image is reconstructed on a grid
of size $m\times m\times m$
from $\mathcal{O}(m^{2})$ projections, each of which contains
values for $\mathcal{O}(m)$ integration spheres, all these
algorithms have complexity of $\mathcal{O}(m^{5})$ operations.
In practical terms, for $m$ of order of a hundred, the
reconstruction time is measured in hours. An example of the
reconstruction in $3D$ using a method based on formula~(\ref
{E:kunyansky}) is shown in Fig.~\ref{F:recon}. Reconstructions
using formulas~(\ref{FPR3da}) or~(\ref{FPR3d}) are quite similar
in terms of stability, accuracy, and computation time.

\subsubsection{Exact filtered backprojection formulas in $2D$}

Exact inversion formulas were obtained for even dimensions
in~\cite{Finch_even}. Denoting by $g$, as before, the spherical
integrals (rather than averages) of $f$, the formulas in $2D$
look as follows:%

\begin{equation}
f(y)=\frac{1}{4\pi^2 R}\Delta\int\limits_{S}\int\limits_
{0}^{2R}g(z,t)\log|t^{2}-|y-z|^{2}|\ dt\ dl(z),
\label{E:Finch2D}%
\end{equation}
or
\begin{equation}
f(y)=\frac{1}{4\pi^2 R}\int\limits_{S}\int\limits_{0}^
{2R}%
\frac{\partial}{\partial t}\left(  t\frac{\partial}{\partial t}
\frac
{g(z,t)}{t}\right)  \log|t^{2}-|y-z|^{2}|\ dt\ dl(z), \label
{E:Finch2Da}%
\end{equation}
where $B$ is a disk of radius $R$ centered at the origin, and $S=
\partial B$ is its boundary.

Another $2D$ inversion formula~\cite{Kunyansky} takes the
following form (again, corrected for a sign):
\begin{equation}
f(y)=-\frac{1}{8\pi}\mathrm{div}\int\limits_{S}n%
(z)h(z,|y-z|)dl(z), \label{E:backpro}%
\end{equation}
where
\begin{align}
h(z,t)  &  =\int\limits_{\mathbb{R}^{+}}\left[  Y_{0}(\lambda t)
\left(
\int\limits_{0}^{2R}J_{0}(\lambda t^{\prime})g(z,t^{\prime})dt^
{\prime
}\right)  \right. \nonumber\\
&  -\left.  J_{0}(\lambda t)\left(  \int\limits_{0}^{2R}Y_{0}
(\lambda
t^{\prime})g(z,t^{\prime})dt^{\prime}\right)  \right]  \lambda d
\lambda,
\label{E:filtration}%
\end{align}
and $J_{0}(t)$ and $Y_{0}(t)$ are the Bessel and Neumann
functions of order $0$. By analyzing the large argument
asymptotics of these functions one
can see~\cite{Kunyansky} that the filtration operator given by
equation~(\ref{E:filtration}) is an analog of the Hilbert
transform.

This reconstruction procedure can be re-written in a
form similar to~(\ref{E:Finch2D}) or~(\ref{E:Finch2Da}).
Indeed, by slightly modifying the original derivation of~(\ref
{E:backpro}),~(\ref{E:filtration}), one can obtain a formula
that would reconstruct a smoothed version $\hat{f}(x,\nu)$ of $f
(x)$ defined by the formula
\[
\hat{f}(x,\nu)=\mathcal{F}^{-1}\left(  |\xi|^{-\nu}\mathcal{F}f
\right) ,\quad0<\nu<1,
\]
where $\mathcal{F,F}^{-1}$ are correspondingly the 2D Fourier
and inverse Fourier transforms. The restriction of $\hat{f}(x,
\nu)$ to the interior of the disk $B$ is recovered by the formula
\begin{equation}
\hat{f}(y,\nu)=-\frac{1}{8\pi}\mathrm{div}\int\limits_{S}
n(z)h_{\nu}(z,|y-z|)dl(z),
\end{equation}
where
\begin{equation}
h_{\nu}(z,t)=\int\limits_{\mathbb{R}^{+}}Y_{0}(\lambda t)\left(
\int\limits_{0}^{2R}J_{0}(\lambda t^{\prime})g(z,t^{\prime})dt^
{\prime
}\right)  -J_{0}(\lambda t)\left(  \int\limits_{0}^{2R}Y_{0}
(\lambda
t^{\prime})g(z,t^{\prime})dt^{\prime}\right)  \lambda^{-\nu}d
\lambda.
\label{newfilt}%
\end{equation}
For $0<\nu<1$, one can change the order of integration in~(\ref
{newfilt}) to obtain
\begin{align}
h_{\nu}(z,t)  &  =\int\limits_{0}^{2R}g(z,t^{\prime})K_{\nu}
(z,t,t^{\prime
})dt^{\prime},\label{inner}\\
K_{\nu}(z,t,t^{\prime})  &  =\int\limits_{\mathbb{R}^{+}}Y_{0}
(\lambda
t)J_{0}(\lambda t^{\prime})\lambda^{-\nu}d\lambda-\int\limits_
{\mathbb{R}^{+}%
}J_{0}(\lambda t)Y_{0}(\lambda t^{\prime})\lambda^{-\nu}d\lambda.
\label{kernel}%
\end{align}
Using~\cite[formula 4.5, p. 211]{Ober}, the integral
$\int\limits_{\mathbb{R}^{+}}Y_{0}(\lambda t)J_{0}(\lambda t^
{\prime}%
)\lambda^{-\nu}d\lambda$ can be integrated exactly, yielding
\[
\int\limits_{\mathbb{R}^{+}}Y_{0}(\lambda t)J_{0}(\lambda t^
{\prime}%
)\lambda^{-\nu}d\lambda=\left\{
\begin{array}
[c]{cc}%
\frac{2^{1-\nu}}{\pi}\Gamma(1-\nu)\frac{t^{-\nu}\cos(\pi\nu)}{|t^
{2}%
-t^{\prime2}|^{1-\nu}}, & t>t^{\prime}\\
-\frac{2^{1-\nu}}{\pi}\Gamma(1-\nu)\frac{t^{-\nu}}{|t^{2}-t^
{\prime2}|^{1-\nu
}}, & t<t^{\prime}%
\end{array}
\right.  .
\]
The expression for the second integral in~(\ref{kernel}) is
derived by interchanging $t$ and $t^{\prime},$ which results in
the formula
\[
K_{\nu}(z,t,t^{\prime})=\left\{
\begin{array}
[c]{cc}%
\frac{2^{1-\nu}}{\pi}\Gamma(1-\nu)\frac{t^{-\nu}\cos(\pi\nu)+(t^
{\prime
})^{-\nu}}{|t^{2}-t^{\prime2}|^{1-\nu}}, & t>t^{\prime}\\
-\frac{2^{1-\nu}}{\pi}\Gamma(1-\nu)\frac{(t^{\prime})^{-\nu}\cos
(\pi
\nu)+t^{-\nu}}{|t^{2}-t^{\prime2}|^{1-\nu}}, & t<t^{\prime}%
\end{array}
\right.
\]
Finally, we substitute the above expression for $K_{\nu}(z,t,t^
{\prime})$ into~(\ref{inner}) and take the limit $\nu\rightarrow0
$, to arrive at the following formulas%
\begin{align*}
f(y)  &  =\frac{1}{2\pi^{2}}\mathrm{div}\int\limits_{S}%
n(z)h_{0}(z,|y-z|)dl(z),\\
h_{0}(z,t)  &  =\int\limits_{0}^{2R}g(z,t^{\prime})\frac{1}
{ { t^{\prime} }^2- t^2  }
dt^{\prime}%
\end{align*}
or
\begin{equation}
f(y)=\frac{1}{2\pi^{2}}\mathrm{div}\int\limits_{S}n%
(z)\left[  \int\limits_{0}^{2R}g(z,t^{\prime})
\frac{1}{  {t^{\prime}}^2 - |y-z|^{2}      }
dt^{\prime}\right]  dl(z). \label{Kunya2D}%
\end{equation}

Similarly to the one appearing in~(\ref{E:Finch2D}) and~(\ref
{E:Finch2Da}), the filtration operator in~(\ref{Kunya2D}) also
involves kernel $\frac{1}%
{ {t^{\prime}}^2-t^2 }.$ If desired, it can be re-written in the
form of a convolution, either by a change of variables $t^{2}
\rightarrow t$, or by noticing that
\[
\frac{2}{ {t^{\prime}}^2 - t^{2} }=\frac{1/t^{\prime}}{t+t^{\prime}}-\frac
{1/t^{\prime}}{t-t^{\prime}}.
\]
This is important from the computational point of view, since it
allows the reduction of the inner integral in~(\ref{Kunya2D})
to the sum of two Hilbert transforms, computational algorithms
for which are well known.

All inversion formulas presented in this section require $
\mathcal{O}(m^{3})$ operations to reconstruct an image on a grid
of size $m\times m$ from $\mathcal{O}(m)$ projections, each
consisting of  $\mathcal{O}(m)$ values
of circular integrals. This coincides with the operation count
required by a classical (non-accelerated) filtered
backprojection algorithm in $2D$.

It is not yet known currently whether formula~(\ref{Kunya2D}) is
equivalent to~(\ref{E:Finch2D}) and~(\ref{E:Finch2Da}). However,
as shown in the previous section, this is not the case for the $3D$
versions of these formulas, and thus this seems unlikely
in the two-dimensional case as well.

Finally, similarly to the filtered backprojection formulas for
the classical $2D$ Radon transform,
the inversion formulas~(\ref{E:Finch2D}),~(\ref{E:Finch2Da}),
and~(\ref{Kunya2D}) are not local. In other words, in order to
recover the value of $f(x)$ for a fixed point $x$, all the values
of $g(z,t)$ have to be known.

\subsection{Series solutions for arbitrary geometries}

%%%%%%%%%%%%%%%%%%%%%%%%%%%%%

Explicit inversion formulas for closed surfaces $S$ different
from spheres have not yet been found\footnote{Planar and cylindrical
observation surfaces, for which such formulas are known
\cite{Natt_new,And,Faw,XFW,XXW}, are not closed.}, except the result of~\cite{AKu} described in the next Section.
There is, however, a
different approach~\cite{Kun_series} that theoretically works
for any closed $S$ and that is practically useful when the
surface is the boundary of a region, in which the spectrum and
eigenfunctions of the Dirichlet Laplacian are known (or could be
effectively approximated numerically).

Let $\lambda_{k}^{2}$ (where $\lambda_k>0$) and $u_{k}(x)$ be
the eigenvalues and normalized eigenfunctions of the Dirichlet
Laplacian $-\Delta_D$ on the interior $\Omega$
of the observation surface $S$:%

\begin{align}
\Delta u_{k}(x)+\lambda_{k}^{2}u_{k}(x)  &  =0,\qquad x\in\Omega,\quad
\Omega\subseteq\mathbb{R}^{n},\label{Helmeq}\\
u_{k}(x)  &  =0,\qquad x\in S =\partial\Omega ,\nonumber\\
||u_{k}||_{2}^{2}  &  \equiv\int\limits_{\Omega}|u_{k}(x)|^{2}dx=1.\nonumber
\end{align}
As before, we would like to reconstruct a compactly supported
function $f(x)$ from the known values of its spherical integrals
$g(z,r)$~(\ref{E:spher_int}).

According to~\cite{Kun_series}, if $f(x)$ is represented as the
sum of the Fourier series
\begin{equation}
f(x)=\sum_{m=0}^{\infty}\alpha_{k}u_{k}(x), \label{fourierser}%
\end{equation}
the Fourier coefficients $\alpha_{k}$ can be reconstructed as
follows:
\begin{equation}
\alpha_{k}=\int_{\partial\Omega}I(z,\lambda_{k})\frac{\partial}
{\partial
n}u_{k}(z)dA(z) \label{serkoef1a}%
\end{equation}
where
\begin{equation}
I(z,\lambda_{k})=\int\limits_{\mathbb{R}^{+}}g(z,r)\Phi_{\lambda_{k}%
}(r)dr,\nonumber
\end{equation}
and $\Phi_{\lambda_{k}}(|x-z|)$ is a free-space rotationally
invariant Green's function of the Helmholtz equation~(\ref
{Helmeq}).

Formula~(\ref{serkoef1a}) is obtained by substituting the
Helmholtz representation for $u_{k}(x)$
\begin{equation}
u_{k}(x)=\int_{\partial\Omega}\Phi_{\lambda_{k}}(|x-z|)\frac
{\partial
}{\partial n}u_{k}(z)ds(z)\qquad x\in\Omega,
\label{helmdiscr}%
\end{equation}
into the expression for the projections $g(z,t).$

This eigenfunction expansion approach requires the knowledge of
the spectrum and eigenfunctions of the Dirichlet Laplacian,
which is available only for some simple domains. However, when
this
information is available, the method yields reliable, robust, and,
in some cases, fast reconstruction. For example,
as it was shown in~\cite{Kun_series}, for the cubic observation
surface $S$, one can compute reconstructions thousands times
faster than by methods based on
explicit inversion formulas of backprojection type discussed
above. The operation count for such an
algorithm is $\mathcal{O}(m^{3}\log m)$, as compared to $\mathcal
{O}(m^{5})$ for the explicit inversion formulas.

\begin{figure}[th]
\begin{center}
%\scalebox{0.7}
%
%
%
{\includegraphics{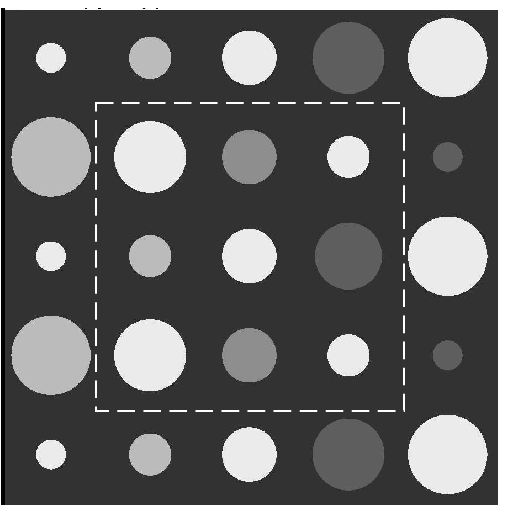} \phantom{aaa} \includegraphics
{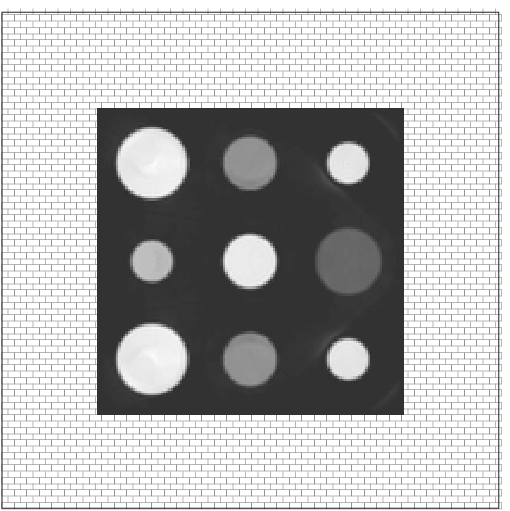}}
\end{center}
\caption{The phantom shown on the left includes several balls located
outside the square acquisition surface $S$, which does not influence
the reconstruction inside $S$ (right).}%
\label{F:good}%
\end{figure}

Another advantage of the series technique is its ability to
''tune out'' the signal coming from outside of $S$. In other
words, unlike the explicit inversion formulas discussed in the
previous sections, the present method enables one to reconstruct
the values of $f(x)$ for all $x$ lying inside
$S$ even in the presence of the sources outside. We illustrate
this property by the reconstruction shown in Fig.~\ref{F:good}.
(The dashed line in the left figure represents surface $S$,
i.e., the location of the detectors.)

\section[Reconstruction: variable speed]{Reconstruction in the case of variable sound speed.}

In this section we consider a more general case of the variable sound speed $v_s(x)$.
Our analysis is valid under previously imposed conditions on this speed,
namely, that $v_s(x)$ is sufficiently smooth, strictly positive, non-trapping,
and $v_s(x)-1$ is compactly supported.

Consider the Hilbert space $H=L^2(\Omega,{v_s}^{-2}(x)dx)$, i.e., the
weighted $L^2$ space with the weight ${v_s}^{-2}(x)$. In this space, the
naturally defined operator
$$
A=-{v_s}^2(x)\Delta
$$
in $\Omega$ with zero Dirichlet conditions on $S$ is self-
adjoint, positive, and has discrete spectrum $\{\lambda_k^2\}
(\lambda_k>0) $ with eigenfunctions $\psi_k(x)\in H$.

We also denote by $E$ the operator of harmonic extension of
functions from $S$ to $\Omega$. I.e., for a function $\phi$ on
$S$ the function $E\phi$ is harmonic inside $\Omega$ and coincides
with $\phi$ on $S$.

Since we are dealing with the unobstructed wave
propagation in the whole space (the surface $S$ is not truly a boundary, but just an
observation surface), and since we assumed that the sound speed
is non-trapping and constant at infinity, the local energy decay
type estimates of~\cite{Vainb,Vainb2} (see also~\cite[Theorem
2.104]{Egorov}) apply. They also lead to the
following reconstruction procedures:

\begin{theorem}\cite{AKu}\label{T:main_speed}
\begin{enumerate}
\item The function $f(x)$ in~(\ref
{E:wave_data}) can be reconstructed inside $\Omega$ as follows:
\begin{equation}\label{E:reconstruction_variable}
f(x)=(Eg|_{t=0})-\int\limits_0^\infty
A^{-\frac12} \sin{(\tau
A^{\frac12})}E(g_{tt})(x,\tau)d\tau.
\end{equation}

\item Function $f(x)$ can be reconstructed inside
$\Omega$ from the data $g$ in~(\ref{E:wave_data}), as the following
$L^2(\Omega)$-convergent series:
\begin{equation}\label{E:coef_variable}
f(x)=\sum\limits_k f_k \psi_k(x),
\end{equation}
where the Fourier coefficients $f_k$ can be recovered using one
of the following formulas:
\begin{equation}\label{E:coef_variable2}
\begin{cases}
f_k=\lambda_k^{-2}g_k(0)-\lambda_k^{-3}\int\limits_0^\infty
\sin{(\lambda_k t)} g_k^{\prime\prime}(t)dt,\\
f_k=\lambda_k^{-2}g_k(0)+\lambda_k^{-2}\int\limits_0^\infty
\cos{(\lambda_k t)} g_k^{\prime}(t)dt, \mbox{ or } \\
f_k=-\lambda_k^{-1}\int\limits_0^\infty \sin{(\lambda_k t)}g_k(t)
dt
=-\lambda_k^{-1}\int\limits_0^\infty \int\limits_S\sin{(\lambda_k
t)}g (x,t)\overline{\frac{\partial \psi_k}{\partial n}(x)}dxdt,
\end{cases}
\end{equation}
where
$$
g_k(t)=\int\limits_{S}g(x,t)\overline{\frac{\partial \psi_k}
{\partial n}(x)}dx
$$
and $n$ denotes the external normal to $S$.
\end{enumerate}
\end{theorem}

\begin{remark} The function $E(g_{tt})$ does not belong to the domain of
the operator $A$. The formula (\ref{E:reconstruction_variable}), however,
still makes sense, since the operator $A^{-\frac12} \sin{(\tau
A^{\frac12})}$ is bounded in $L^2$. \end{remark}

This theorem in the particular case of the constant sound speed,
implies the eigenfunction expansion procedure of~\cite
{Kun_series} described in the previous section. However,
unlike~\cite{Kun_series},
it also applies to the variable speed situation
and it does not require knowledge of a whole space Green's
function. Similarly to the method of~\cite{Kun_series}
discussed in the preceding section, this procedure yields
correct reconstruction inside the domain, even if a part of the
source lies outside.

\section{Partial data. ``Visible'' and ``invisible''
singularities}

\label{S:partial}
%%%%%%%%%%%%%%%%%%%%%%%%%%%%%

One can find a more detailed discussion of this issue for
TAT in~\cite{KuKu,XWAK}. Here we provide only a brief summary.

Uniqueness of reconstruction does not necessarily mean the
possibility of practical reconstruction, since the
reconstruction procedure can sometimes be unstable. This
is true, for instance,  in problems of electrical impedance
tomography, and in incomplete data problems of X-ray tomography
and TAT~\cite{Kuch_AMS05,KuchQuinto,Natt_old,Natt_new}.

Microlocal analysis done in~\cite{LQ,Quinto} (see also~\cite
{Pal_book}) shows which
parts of the wave front of a function $f$ can be recovered
from its partial $X$-ray or TAT data (see also~\cite{XWAK}
for a practical discussion). We describe this result
in an imprecise form (see~\cite{LQ} for precise
formulation), restricted to the case of jump singularities
(tissue interfaces) only.

\begin{figure}[th]
\begin{center}
%\scalebox{0.7}
\includegraphics{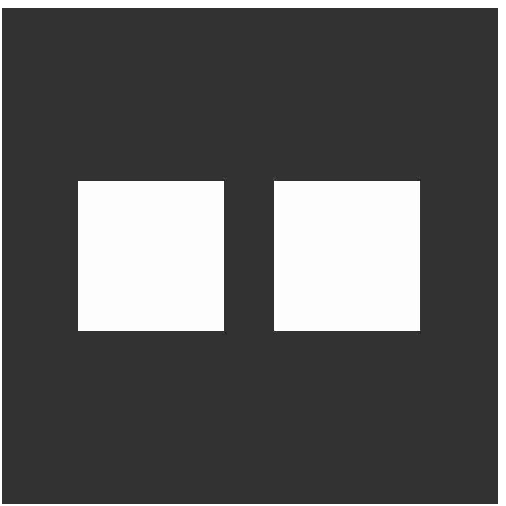}
%\scalebox{0.7}
\phantom{aaa}
\includegraphics{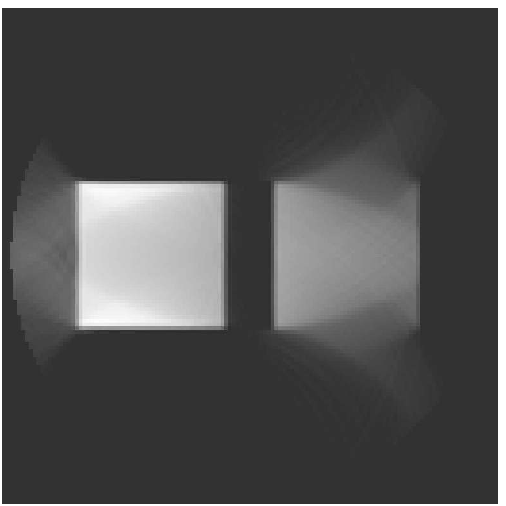}
\end{center}
\caption{Effect of incomplete data: the phantom (left) and its incomplete data
reconstruction.}
\label{F:incomplete}%
\end{figure}

According to~\cite{LQ,Quinto}, for such singularities
a part of the interface is stably recoverable
(dubbed ``visible'' or ``audible''), if for each
point of the interface there exists a sphere
centered at $S$ and tangent to the interface at this
point. Otherwise, the interface will be blurred away (even
if there is a uniqueness of reconstruction theorem).
Indeed, if all spheres of integration are transversal to
the interface, the integration smooths the singularity,
and thus reconstruction of this interface becomes unstable.
The Figure~\ref{F:incomplete} shows an example of a
reconstruction from incomplete spherical mean data.
The simulated transducers in this experiment were
located along a $180^{o}$ circular arc (the left half of
a large circle surrounding the squares). In this figure
the sides of the squares that are not touched
tangentially by circles centered on $S$ are noticeably
blurred; any kind of de-blurring technique would not
be stable in this context.

%%%%%%%%%%%%%%%%%%%%%
\section{Range conditions}
%%%%%%%%%%%%%%%%%%%%%

This paper would not be complete without mentioning
the intimate relationship of inversion problems with
range conditions.
Indeed, as it has already been mentioned, recovery of $f$
from the data $g$ is impossible, if considered as an inverse
problem for the wave
equation problem inside the cylinder $S\times\RR^+$.
The possibility of inversion depends upon the fact that
the solution of the wave equation lives in the whole
space, and $S$ is just the observation surface, rather
than a true boundary. In other words, the data $g(x,t)$
comes from a very small (infinite co-dimension) subspace
in any natural function space on the lateral boundary $S
\times\RR^+$. Thus, range conditions must play a
significant role. Indeed, they lead the authors of~\cite
{AKu} to their results. We thus provide here a brief
sketch of range results, following essentially the
corresponding section of~\cite{KuKu}.

As it has just been mentioned, the ranges of Radon type transforms,
including the spherical mean operator, are usually of infinite co-
dimension in natural function spaces (in other words, ideal data
should satisfy infinitely many consistency conditions). Information about
the range is important for many theoretical and practical
purposes (reconstruction algorithms, error
corrections, incomplete data completion, etc.), and has
attracted a lot of attention (e.g.,
\cite{Leon_Radon,GGG1,GGG,GelfVil,Helg_Radon,Helg_groups,Kuch_AMS05,KucLvin,KucLvin2,KuchQuinto,Lvin,Natt_old,Natt_new,Novikov,Pal_book,Ponomarev}).

For example, functions $g$ from the range of the standard Radon
transform
\[
f(x)\to g(s,\omega)=\int\limits_{x\cdot\omega=s}f(x)dx, |\omega|
=1,
\]
satisfy two types of conditions:

\begin{enumerate}
\item \emph{evenness}: $g(-s,-\omega)=g(s,\omega)$

\item \emph{moment conditions}: for any integer $k\geq0$, the $k
$th moment
\[
G_{k}(\omega)=\int\limits_{-\infty}^{\infty} s^{k} g(\omega,s)ds
\]
extends from the unit circle of vectors $\omega$ to a
homogeneous polynomial of degree $k$ in $\omega$.
\end{enumerate}

Although for the Radon transform the evenness condition
seems to be ``trivial'', while the moment conditions seem
to be the most important, this perception is misleading.
Indeed, for more general transforms of Radon type it is
often easier to find analogs of the moment conditions,
while counterparts of the evenness conditions could be
elusive (see~\cite{Kuch_AMS05,KucLvin,KucLvin2,Natt_old,Natt_new,Novikov}).
This is exactly what happens with the spherical
mean transform $R_{S}$.

An analog of the moment conditions was first present
implicitly in~\cite{AQ,LP1,LP2} and explicitly
formulated as such in~\cite{etti,Patch}:

\textbf{Moment conditions} \emph{on data $g(x,r)=R_{S} f
(x,r)$ in $\RR^n$ are: for any
integer $k\geq0$, the moment
\[
M_{k}(x)=\int\limits_{0}^{\infty} r^{2k+n-1} g(x,r)dr, x\in S
\]
can be extended from $S$ to a (non-homogeneous) polynomial
$Q_{k}(x)$ of
degree at most $2k$.}

These conditions are incomplete, and infinitely many
others, which play the role of an analog of evenness, need
to be added.

Complete range description for $R_{S}$ when $S$ is a
sphere in $2D$ was found in~\cite{AmbKuc_rang} and then in
odd dimensions in~\cite{FR}. They were then extended to
any dimension and provided several interpretations
in~\cite{AKQ}. These conditions, which happen to be
intimately related to PDEs and spectral theory, are
described below.

\begin{figure}[th]
\begin{center}
%\scalebox{0.7}
{\includegraphics{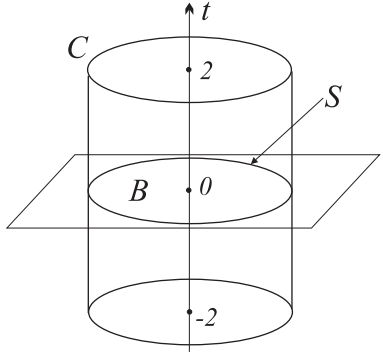}}
\end{center}
\caption{An illustration to the range description.}%
\label{F:}%
\end{figure}

Let $B$ be the unit ball in $\mathbb{R}^{n}$, $S=\partial B$
the unit sphere, and $C$ the cylinder $B\times[0,2]$
(see Fig.~\ref{F:}).

Consider the spherical mean operator $R_{S}$:
\[
R_{S}f(x,t)=G(x,t)=\int_{|y|=1}f(x+ty)dA(y).
\]
If $G(x,t)$ is defined by the same formula for all $x\in
\mathbb{R}^{n}$, then it satisfies Darboux (Euler-Poisson-
Darboux) equation~\cite{Asgeirsson,CH,John}
\[
G_{tt}+(n-1)t^{-1}G_{t}=\Delta_{x} G.
\]
Inside the cylinder $C$, $G(x,t)$ vanishes when $t\geq2$
(since the spheres of integration do not intersect the
support of the function when $t\geq2$).

\begin{theorem}
\cite{AKQ}\label{T:AKQ} The following four statements
are equivalent for any function $g\in
C^{\infty}_{0} (S\times[0,2])$, where $S$ is a sphere:

\begin{enumerate}
\item  Function $g$ is representable as $R_{S} f$ for some $f\in
C^{\infty
}_{0}(B)$.

\item
\begin{enumerate}
\item  The moment conditions are satisfied.

\item  The solution $G(x,t)$ of the interior Darboux problem
satisfies the
condition
\[
\lim\limits_{t \to0}\int\limits_{B} \frac{\partial G}{\partial
t}%
(x,t)\phi(x)dx=0
\]
for any eigenfunction $\phi(x)$ of the Dirichlet Laplacian in $B
$.
\end{enumerate}

\item
\begin{enumerate}
\item  The moment conditions are satisfied.

\item  Let $-\lambda^{2}$ be an eigenvalue of Dirichlet
Laplacian in $B$ and
$\psi_{\lambda}$ the corresponding eigenfunction. Then the
following
orthogonality condition is satisfied:
\begin{equation}
\int\limits_{S\times[0,2]} g(x,t)\partial_{\nu}\psi_{\lambda}(x)
j_{n/2-1}%
(\lambda t)t^{n-1}dxdt=0.
\end{equation}
Here $j_{p}(z)=c_{p}\frac{J_{p}(z)}{z^{p}}$ is the so called
spherical Bessel function.
\end{enumerate}

\item
\begin{enumerate}
\item  The moment conditions are satisfied.

\item  Let $\widehat{g}(x,\lambda)=\int g(x,t)j_{n/2-1}(\lambda
t)t^{n-1}dt$.
Then, for any $m\in\mathbb{Z}$, the $m^{th}$ spherical harmonic
term
$\widehat{g}_{m}(x,\lambda)$ of $\widehat{g}(x,\lambda)$
vanishes at non-zero
zeros of Bessel function $J_{m+n/2-1}(\lambda)$.
\end{enumerate}
\end{enumerate}
\end{theorem}

One can make several important comments concerning this result (see~\cite{AKQ}
for a detailed discussion). In all of the remarks below, except the
third one, the observation surface $S$ is assumed to be a sphere.

\begin{enumerate}
\item  If the dimension $n$ is odd, then conditions (b) alone suffice
for the complete range description,
and thus they imply the moment conditions as well.
(A similar earlier result was established for a related transform
in~\cite{FR}.) It is not clear at the moment whether this is holds
true in even dimensions.

\item  The range descriptions for $R_S$ work in Sobolev scale, i.e.
they describe the range of the operator $R_S: H^{s}_{comp}(B)\mapsto
H^{s+(n-1)/2}_{comp}(S\times\RR^+)$. (This uses a recent work by
Palamodov~\cite {Palam_funk}). Notice that in this result it is
assumed that the function $f$ vanishes in a neighborhood of $S$,
while in the previous theorem it was allowed for the support of $f$
to reach all the way to the sphere $S$.

\item  If $S$ is not a sphere, but the boundary of a bounded domain,
the range conditions 2 and 3 of the previous Theorem
are still necessary for the data $g$ to belong to the range of $R_S$.
They, however, might no longer suffice for $g$ to belong to the range.

\item A different wave equation approach to the range descriptions can
be found in~\cite{FR}.

\end{enumerate}

%%%%%%%%%%%%%%%%%%%%%
\section{Concluding remarks}
%%%%%%%%%%%%%%%%%%%%%

\subsection{Uniqueness}

As it has already been mentioned, the uniqueness questions relevant for TAT
applications are essentially resolved. However, the
mathematical understanding of the uniqueness problem
for the restricted spherical mean operators $R_{S}$ is
still unsatisfactory and open problems abound~\cite{AQ,KuKu}.
For instance, very little is known for the case of
functions without compact support. The
main known result is of~\cite{ABK}, which describes for
which values of $1\leq
p \leq\infty$ the uniqueness result still holds:
\begin{theorem}
\cite{ABK} Let $S$ be the boundary of a bounded domain in $
\mathbb{R}^{n}$ and
$f\in L^{p}(\mathbb{R}^{n})$ such that $R_{S} f\equiv0$. If $p
\leq 2n/(n-1)$,
then $f\equiv0$ (and thus $S$ is injectivity set for this
space). This fails
for any $p>2n/(n-1)$.
\end{theorem}

The three- and higher-dimensional uniqueness problem for
non-closed observation surface $S$ is also still open~\cite{AQ,KuKu}.

\subsection{Inversion}

Albeit closed form (backprojection type) inversion formulas are
available now
for the cases of $S$ being a plane (and object on one side from
it), cylinder,
and a sphere, there is still some mystery surrounding this
issue. For instance, it would be interesting to understand whether
(closed form, rather than series expansion) backprojection type
inversion formulas could be written for non-spherical observation
surfaces $S$ and/or in the presence of a non-uniform
background $v_s(x)$. The results presented in Section 1.5 seem to
be the first step in this direction.

The I.~Gelfand's school of integral geometry has
developed a powerful technique of the so called $\kappa$
operator, which provides a general approach to inversion
and range descriptions for transforms of Radon type~\cite{GGG1,GGG}.
 In particular, it has been applied to
the case of integration over various collections
(``complexes'') of spheres in~\cite{GGG,Gi}. This
consideration seems to suggest that one should not expect
explicit closed form inversion formulas for $R_{S}$ when $S
$ is a sphere. However, such formulas were discovered
in~\cite{FPR,Finch_even,Kunyansky}.
This apparent controversy
(still short of contradiction) has not been resolved completely
yet.

B.~Rubin has recently discovered an alternative interesting approach to
inversion formulas of the type of (\ref{FPR3da})-(\ref{FPR3d}) for the case
when $S$ is a sphere. It relies upon the idea of regarding the spherical
mean operator as a member of a broader family of operators \cite{Rubin}.

In $3D$, if the sound speed is constant, the Huygens' principle applies,
i.e. the pressure $p(t,x)$ inside $S$ becomes equal to zero for any
time $T$ larger than the time required for sound to cross the domain.
Thus, imposing zero conditions on $p(t,x)$ and $p_t(t,x)$ at $t=T$ and
solving the wave equation (\ref{E:wave_data}) back in time with the
measured data $g$ as the boundary values, one recovers at $t=0$ the source
$f(x)$. This method has been implemented in \cite{Burgh}. Although in even
dimensions or in presence of sound speed variations, Huygens'
principle does not apply, one can find good approximate solutions by a similar
approach \cite{BGK}.

A different approach to TAT inversion is suggested in \cite{Klibanov}. It
is based on using not only the measured data $g$ on $S\times\RR^+$, but
also the normal derivative of the pressure $p$ on $S$. Since this normal
derivative is not measured, finding it would require solving the exterior
problem first and deriving the normal derivative from there. Feasibility
and competitiveness of such a method for TAT is not clear at the moment.
\subsection{Stability}

Stability of inversion when $S$ is a sphere surrounding
the support of $f(x)$ is the same as for the standard
Radon transform, as the results of~\cite{AKQ,KuKu,Palam_funk}
show. However, if the support reaches outside, in spite of
Theorem~\ref{T:uniqueness} that
guarantees uniqueness of reconstruction, stability for some
parts of $f(x)$ lying outside $S$ does not hold anymore.
See~\cite{AKQ,KuKu,LQ,Quinto} for details.

\subsection{Range}

The range conditions 2 and 3 of Theorem~\ref{T:AKQ}
are necessary also for non-spherical closed surfaces
$S$ and for functions with support outside $S$. They, however,
are not expected to be sufficient, since the arising
instabilities indicate that one might expect non-closed ranges
in some cases.

%%%%%%%%%%%%%%%%%%%%%
%
%
%
%
%
%

\section*{Acknowledgments}

%%%%%%%%%%%%%%%%%%%%%
%
%

The work of the second author was partially supported by the NSF
DMS grants 0604778 and 0648786. The third author was partially supported by
the DOE grant DE-FG02-03ER25577 and NSF DMS grant 0312292. The
work was partly done when the first two authors were visiting
the Isaac Newton Institute for Mathematical Sciences (INI) in
Cambridge. The
authors express their gratitude to the NSF, DOE, and INI for this
support. They also thank G.~Ambartsoumian, G.~Beylkin, D.~Finch, A.~Greenleaf,
M.~Klibanov, V.~Palamodov, P.~Stefanov, B.~Vainberg, and E.~Zuazua for information,
and the reviewers and the editor for useful comments.

%%%%%%%%%%%%%%%%%%%%%

%\end{doublespace}

\begin{thebibliography}{999}

\bibitem{Kruger} Kruger, R.~A., Liu, P., Fang, Y.~R. and
Appledorn, C.~R. 1995. Photoacoustic ultrasound (PAUS) reconstruction
tomography.
\textit{Med. Phys.}
22:1605-09.

\bibitem{Kruger03} Kruger, R. A., Kiser, W. L., Reinecke, D. R. and
Kruger, G. A. 2003.
Thermoacoustic computed tomography using a conventional linear transducer array.
\textit{Med. Phys.}
30(5):856-60.

\bibitem{Oraev96} Oraevsky, A.~A., Esenaliev, R.~O., Jacques,  and S.~L.
Tittel, F.~K. 1996.
Laser optoacoustic tomography for medical diagnostics principles.
\textit{Proc. SPIE} 2676:22.

\bibitem{Oraev} Oraevsky, A.~A.  and  Karabutov, A.~A. 2002.
In \textit{Handbook of Optical Biomedical Diagonstics},
edited by V. V. Tuchin, SPIE, Bellingham, WA, Chap. 10.

\bibitem{Oraev2}Oraevsky A.~A. and  Karabutov, A.~A. 2003
Optoacoustic Tomography, Ch. 34
In \textit{Biomedical Photonics Handbook}, edited by
T. Vo-Dinh, CRC, Boca Raton, FL, Chap. 34, {\bf 34}-1 -- {\bf
34}-34.

\bibitem{PAT} Wang, X., Pang, Y., Ku, G. et al. 2003.
Noninvasive
laser-induced photoacoustic tomography for structural and
functional \emph{in vivo} imaging of the brain.
\textit{Nature Biotechnology}
21(7):803-806.


\bibitem{Wang_book} Wang, L.~V. and Wu, H. 2007.
\textit{Biomedical Optics. Principles and Imaging}.
Wiley-Interscience.

\bibitem{MXW_review} Xu, M. and Wang, L.-H.~V. 2006.
Photoacoustic imaging in biomedicine.
\textit{Review of Scientific Instruments} 77:041101-01 -- 041101-22.

\bibitem{KuKu} Kuchment, P. and Kunyansky, L. 2007. Mathematics of
thermoacoustic and photoacoustic tomography. Preprint
arXiv:0704.0286v1 [math.AP], submitted.

\bibitem{LQ} Louis, A.~K. and Quinto, E.~T. 2000. Local tomographic
methods in Sonar. In \textit{Surveys on solution methods for
inverse problems}, 147-154,
Vienna:Springer.

\bibitem{NC} Nolan C.~J. and Cheney, M. 2002.
Synthetic aperture inversion.
\textit{Inverse Problems} 18:221-235.

\bibitem{Beylkin} Beylkin, G. 1984
The inversion problem and applications of the
generalized Radon transform.
\textit{Comm. Pure Appl. Math.} 37:579-599.

\bibitem{MXW1} Xu, M., and  Wang, L.-H.~V. 2002.
Time-domain reconstruction for
thermoacoustic tomography in a spherical geometry.
\textit{IEEE Trans. Med. Imag.} 21:814-822.

\bibitem{AKu} Agranovsky,  M. and  Kuchment, P. 2007. Uniqueness of
reconstruction and an inversion procedure
for thermoacoustic and photoacoustic tomography.
Preprint, arXiv:0706.0598.

\bibitem{JinWang} Jin, X. and  Wang, L.~V. 2006.
Thermoacoustic tomography with correction for acoustic speed variations.
\textit{Physics in Medicine and Biology}
51:6437-48.

\bibitem{haltmaier_large} Haltmeier, M., Burgholzer, P.,
Paltauf, G. and  Scherzer, O. 2004.
Thermoacoustic computed tomography with large planar
receivers.
\textit{Inverse Problems} 20:1663-73.

\bibitem{haltmaier_fabri} Burgholzer, P., Hofer, C.,  Matt, G.~J. et al. 2006.
Thermoacoustic tomography
using a fiber-based Fabry-Perot interferometer as an integrating line
detector.
\textit{Proc. SPIE} 6086:434-442.


\bibitem{And} Andersson, L.-E. 1988.
On the determination of a function from spherical averages.
\textit{SIAM J. Math. Anal.}
19(1):214-32.

\bibitem{Faw} Fawcett, J. A. 1985.
Inversion of $n$-dimensional spherical averages.
\textit{SIAM J. Appl. Math.}
45(2):336-41.

\bibitem{Natt_new} Natterer, F. and W\"{u}bbeling, F. 2001.
\textit{Mathematical Methods in Image Reconstruction}.
Monographs on Mathematical Modeling and Computation v. 5.
Philadelphia: SIAM.


\bibitem{XFW} Xu, Y., Feng, D. and  Wang, L.-H. V. 2002.
Exact frequency-domain reconstruction for thermoacoustic
tomography: I. Planar geometry.
\textit{IEEE Trans. Med. Imag.}
21:823-28.

\bibitem{XXW}Xu, Y., Xu, M.  and  Wang, L.-H.  V. 2002.
Exact frequency-domain reconstruction for thermoacoustic
tomography: II. Cylindrical geometry.
\textit{IEEE Trans. Med. Imag.}
21:829-33.

\bibitem{Diebold} Diebold, G.~J., Sun, T., Khan, M.~I. 1991.
Photoacoustic monopole radiation in one, two, and three dimensions.
\textit{Phys. Rev. Lett.}
67(24):3384-87.

\bibitem{Tam} Tam, A.~C. 1986. Applications of photoacoustic sensing techniques.
\textit{Rev. Mod. Phys.}
58(2):381-431.

\bibitem{ABK} Agranovsky, M.,  Berenstein, C. and  Kuchment, P.  1996.
Approximation by spherical waves in $L^{p}$-spaces.
\textit{J. Geom. Anal.}
6(3):365-83.

\bibitem{AKQ} Agranovsky, M., Kuchment, P., and Quinto, E.~T.    2007.
Range descriptions for the spherical mean Radon transform.
\textit{J. Funct. Anal.}
248:344-86.

\bibitem{AQ} Agranovsky M. and  Quinto, E. T. 1996.
Injectivity sets for the Radon
transform over circles and complete systems of radial functions.
\textit{J. Funct. Anal.}
139:383-414.

\bibitem{AmbKuc_inj} Ambartsoumian G. and Kuchment, P. 2005.
On the injectivity of the circular Radon transform.
\textit{Inverse Problems}
21:473-85.

\bibitem{FPR} Finch, D., Patch, S. and Rakesh. 2004.
Determining a function from its mean values over a family of spheres.
\textit{SIAM J. Math. Anal.}
35(5):1213-40.

\bibitem{CH} Courant R. and  Hilbert, D. 1962.
\textit{Methods of Mathematical Physics, Volume II Partial Differential Equations},
New York: Interscience.




\bibitem{FR2} Finch, D. and Rakesh. 2007.
Recovering a function from its spherical mean values in two and three
dimensions. In this volume.



\bibitem{FR3} Finch, D. and Rakesh. 2007.
The spherical mean value operator
with centers on a sphere. Preprint. To appear in \textit{Inverse Problems}.

\bibitem{Kuc93}Kuchment, P. 1993. Unpublished.

\bibitem{Kuch_AMS05}Kuchment, P.
Generalized Transforms of Radon Type and Their Applications. In
~\cite{OlafQuinto}, 67-91.

\bibitem{Tataru} Tataru, D. 1995.
Unique continuation for solutions to PDEs;
between H\"{o}rmander's theorem and Holmgren's theorem.
\textit{Comm. PDE} 20:814-22.


\bibitem{Norton1} Norton, S. J. 1980.
Reconstruction of a two-dimensional reflecting
medium over a circular domain: exact solution.
\textit{J. Acoust. Soc. Am.}
67:1266-73.

\bibitem{Norton2} Norton, S. J. and Linzer, M. 1981.
Ultrasonic reflectivity imaging in three dimensions:
exact inverse scattering solutions for plane,
cylindrical, and spherical apertures.
\textit{IEEE Trans. on Biomed. Eng.}
28:200-202.

\bibitem{Leon_Radon} Ehrenpreis, L. 2003.
\textit{The Universality of the Radon Transform}
Oxford Univ. Press .

\bibitem{GGG} Gelfand, I.,  Gindikin, S. and Graev, M. 2003.
\textit{Selected Topics in Integral Geometry},
Transl. Math. Monogr. v. 220.
Providence: Amer. Math. Soc.


\bibitem{Helg_Radon} Helgason, S. 1980.
\textit{The Radon Transform}.
Basel: Birkh\"{a}user.

\bibitem{Natt_old} Natterer, F. 1986.
\textit{The mathematics of computerized tomography}.
New York: Wiley.

\bibitem{PopSush} Popov, D. A. and Sushko, D. V. 2002.
A parametrix for the problem of
optical-acoustic tomography.
\textit{Dokl. Math.}
65(1):19-21.

\bibitem{PopSush2} Popov, D. A. and  Sushko, D. V. 2004.
Image restoration in optical-acoustic tomography.
\textit{Problems of Information Transmission}
40(3):254-278.

\bibitem{XWAK} Xu, Y.,  Wang, L., Ambartsoumian, G. and
Kuchment, P.  2004.
Reconstructions in limited view thermoacoustic tomography.
\textit{Medical Physics}
31(4):724-33.

\bibitem{Finch_even}Finch, D. Haltmeier, M. and Rakesh.   2007.
Inversion of spherical means and the wave equation in even dimensions.
Preprint arXiv math.AP/0701426.

\bibitem{Kunyansky} Kunyansky, L. 2007.
Explicit inversion formulae for the spherical mean Radon transform.
\textit{Inverse Problems}
23:373-83.

\bibitem{MXW2} Xu, M. and  Wang, L.-H. V. 2005.
Universal back-projection algorithm for
photoacoustic computed tomography.
\textit{Phys. Rev. E}
71:016706.

\bibitem{AmbPatch} Ambartsoumian, G. and Patch, S. 2006
Thermoacoustic tomography - implementation of exact backprojection formulas.
Preprint, arXiv:math.NA/0510638.

\bibitem{Ober} Oberhettinger, F. 1972.
\textit{Tables of Bessel Transforms}.
New York: Springer.

\bibitem{Kun_series} Kunyansky, L. 2007.
A series solution and a fast algorithm for the inversion of the
spherical mean Radon transform.
Preprint, arXiv math.AP/0701236.  To appear in \textit{Inverse Problems}.

\bibitem{Vainb} Vainberg, B. 1975.
The short-wave asymptotic behavior
of the solutions of stationary problems, and the
asymptotic behavior as $t\to\infty$ of the solutions of
nonstationary problems.
\textit{Russian Math. Surveys}
30(2):1-58.

\bibitem{Vainb2} Vainberg, B. 1989.
\textit{Asymptotics methods in the Equations of Mathematical Physics}.
Gordon and Breach,  (Translation of the Russian 1982 edition).

\bibitem{Egorov} Egorov, Yu. V. and Shubin, M. A. 1992.
Linear Partial Differential Equations. Foundations of the Classical Theory.
In \textit{Partial Differential Equations. I.},
ed. Yu. V. Egorov and M. A. Shubin, Encyclopaedia of
Mathematical Sciences, 30:1-259. Springer Verlag.

\bibitem{KuchQuinto}Kuchment, P. and Quinto, E. T.
Some problems of integral geometry arising in tomography.
Chapter XI in~\cite {Leon_Radon}.

\bibitem{Quinto} Quinto, E. T. 1993.
Singularities of the X-ray transform and limited data tomography in $\mathbb{R}^{2}$
and $\mathbb{R}^{3}$.
\textit{SIAM J. Math. Anal.}
24:1215-25.

\bibitem{Pal_book} Palamodov, V. P. 2004.
\textit{Reconstructive Integral Geometry}.
Basel: Birkh\"{a}user.

\bibitem{GGG1} Gelfand, I., Gindikin, S. and Graev, M. 1980.
Integral geometry in affine and projective spaces.
\textit{J. Sov. Math.}
18:39-167.

\bibitem{GelfVil} Gelfand, I., Graev, M. and Vilenkin, N. 1965.
\textit{Generalized Functions, v. 5: Integral Geometry and Representation Theory}.
Acad. Press.

\bibitem{Helg_groups}Helgason, S. 2000.
\textit{Groups and Geometric Analysis}.
Providence: Amer. Math. Soc.

\bibitem{KucLvin} Kuchment, P. and Lvin, S. 1990.
Paley-Wiener theorem for the exponential Radon transform.
\textit{Acta Applicandae Mathematicae}
18:251-60.

\bibitem{KucLvin2} Kuchment, P. and Lvin, S. 1991.
The Range of the Exponential Radon Transform.
\textit{Soviet Math Dokl.}
42(1):183-184.

\bibitem{Lvin} Lvin, S. 1994.
Data correction and restoration in emission tomography.
In \textit{Tomography, Impedance Imaging, and Integral Geometry},
ed. E.T. Quinto, M. Cheney, and P. Kuchment, 149-155,
Lectures in Appl. Math., vol. 30. Providence: AMS.


\bibitem{Novikov} Novikov, R. 2002.
On the range characterization for the
two-dimensional attenuated X-ray transform.
\textit{Inverse Problems}
18:677-700.

\bibitem{Ponomarev} Ponomarev, I. 1995.
Correction of emission tomography data.
Effects of detector displacement and non-constant sensitivity.
\textit{Inverse Problems}
10:1-8.

\bibitem{LP1} Lin, V. and Pinkus, A. 1993.
Fundamentality of ridge functions.
\textit{J. Approx. Theory}
75:295-311.

\bibitem{LP2}Lin V. and Pinkus, A. 1994.
Approximation of multivariate functions
In \textit{Advances in computational mathematics},
ed. H. P. Dikshit and C. A. Micchelli, 1-9,
World Sci. Publ.

\bibitem{etti} Bouzaglo-Burov, E. 2005.
Inversion of spherical Radon transform, methods and numerical experiments.
MS Thesis,
Bar-Ilan Univ. 1-30. (In Hebrew)

\bibitem{Patch}Patch, S. K. 2004.
Thermoacoustic tomography - consistency conditions and the partial scan problem.
\textit{Phys. Med. Biol.}
49:1-11.

\bibitem{AmbKuc_rang} Ambartsoumian, G. and Kuchment, P. 2006.
A range description for the planar circular Radon transform.
\textit{SIAM J. Math. Anal.}
38(2):681-92.

\bibitem{FR} Finch, D. and Rakesh. 2006.
The range of the spherical mean value operator for functions supported in a ball.
\textit{Inverse Problems}
22:923-38.

\bibitem{Asgeirsson} Asgeirsson, L. 1937.
\"{U}ber eine Mittelwerteigenschaft von
L\"{o}sungen homogener linearer partieller
Differentialgleichungen zweiter
Ordnung mit konstanten Koeffizienten.
\textit{Ann. Math.}
113:321-46.

\bibitem{John} John, F. 1971.
\textit{Plane Waves and Spherical Means, Applied to Partial Differential Equations.}
Dover.

\bibitem{Palam_funk} Palamodov, V.
Remarks on the general Funk-Radon transform
and thermoacoustic tomography.
Preprint, arXiv:math/0701204.


\bibitem{Gi} Gindikin, S. 1995.
Integral geometry on real quadrics.
In \textit{\ Lie groups and Lie algebras: E. B. Dynkin's Seminar},
23-31, Amer. Math. Soc. Transl. Ser. 2, 169,
Providence: Amer. Math. Soc.

\bibitem{Rubin} Rubin, B. 2007. Private communication.

\bibitem{Burgh} Burgholzer, P., Matt, G., Haltmeier, M. \&
Patlauf, G. 2007.
Exact and approximate imaging methods for
photoacoustic tomography using an arbitrary detection surface.
\textit{Phys. Rev. E} 75:046706.

\bibitem{BGK} Bangerth, W., Georgieva-Hristova, Y. \& Kuchment, P. 2007.
On reconstruction in thermoacoustic tomography with variable speed,
in preparation.

\bibitem{Klibanov} Clason, C. and Klibanov, M. 2007.
Quasireversibility method
in thermoacoustic tomography in heterogeneous medium.
Preprint.

\bibitem{OlafQuinto} Olafsson, G. and Quinto, E. T. (Editors), 2006.
\textit {The Radon Transform, Inverse Problems, and Tomography. American
Mathematical Society Short Course January 3-4, 2005, Atlanta,
Georgia}, Proc. Symp. Appl. Math., v. 63, AMS, RI .

\end{thebibliography}
\end{document}